\RequirePackage{etex}

\documentclass[10pt]{amsart}
\usepackage[applemac]{inputenc}

\usepackage[dvipsnames,svgnames,x11names,hyperref]{xcolor}
\usepackage[hyphens]{url}
\usepackage[colorlinks=true,linkcolor=Maroon,citecolor=Maroon,urlcolor=blue,hypertexnames=false,linktocpage]{hyperref}
\usepackage{bookmark}
\usepackage{amsthm,thmtools,amssymb,amsmath,amscd}
\RequirePackage{etex}

\usepackage{fancyhdr}
\usepackage{esint}
\usepackage{enumerate}

\usepackage{pictexwd,dcpic}
\usepackage{graphicx}
\usepackage{caption}
\swapnumbers
\declaretheorem[name=Theorem,numberwithin=section]{thm}
\declaretheorem[name=Remark,style=remark,sibling=thm]{rem}
\declaretheorem[name=Lemma,sibling=thm]{lemma}
\declaretheorem[name=Proposition,sibling=thm]{prop}

\numberwithin{equation}{section}

\usepackage{autonum}

\newcommand{\ti}{\tilde}

\newcommand{\cn}{\colon}
\newcommand{\sub}{\subset}
\newcommand{\ov}{\overline}
\newcommand{\mr}{\mathring}

\newcommand{\bbR}{\mathbb{R}}

\newcommand{\bbS}{\mathbb{S}}

\newcommand{\bbM}{\mathbb{M}}

\newcommand{\8}{\infty}

\newcommand{\al}{\alpha}
\newcommand{\be}{\beta}
\newcommand{\ga}{\gamma}
\newcommand{\de}{\delta}
\newcommand{\ep}{\epsilon}

\newcommand{\si}{\sigma}
\newcommand{\Si}{\Sigma}

\newcommand{\vt}{\vartheta}
\newcommand{\Om}{\Omega}
\newcommand{\De}{\Delta}
\newcommand{\Ga}{\Gamma}

\newcommand{\cH}{\mathcal{H}}

\newcommand{\cU}{\mathcal{U}}

\newcommand{\del}{\partial}
\newcommand{\n}{\nabla}

\newcommand{\ip}[2]{\left\langle #1,#2 \right\rangle}
\newcommand{\fr}[2]{\frac{#1}{#2}}
\newcommand{\tfr}[2]{\tfrac{#1}{#2}}
\newcommand{\x}{\times}

\DeclareMathOperator{\dive}{div}

\DeclareMathOperator{\const}{const}

\DeclareMathOperator{\Rc}{Rc}

\newcommand{\pf}[1]{\begin{proof}#1 \end{proof}}
\newcommand{\eq}[1]{\begin{equation}\begin{alignedat}{2} #1 \end{alignedat}\end{equation}}

\newcommand{\br}[1]{\left(#1\right)}
\newcommand{\abs}[1]{\lvert #1\rvert}

\newcommand{\hp}{\hphantom}
\newcommand{\q}{\quad}

\begin{document}
\title[Stability for Serrin's and Alexandroff's theorem]{Stability for Serrin's problem and Alexandroff's theorem in warped product manifolds}
\author[J. Scheuer]{Julian Scheuer}
\address{Goethe-Universit\"at, Institut f\"ur Mathematik, Robert-Mayer-Str.~10, 60629 Frankfurt, Germany}
\email{scheuer@math.uni-frankfurt.de}
\author[C. Xia]{Chao Xia}
\address{School of Mathematical Sciences, Xiamen University, 361005,
Xiamen, P.R. China}
\email{chaoxia@xmu.edu.cn}

\thanks{CX is  supported by the  NSFC (Grant No. 11871406, 12271449).}
\keywords{Constant mean curvature; Heintze-Karcher inequality; Serrin's problem; Stability; Warped product manifold.}
\date{\today}
\begin{abstract}
We prove quantitative versions for several results from geometric partial differential equations. Firstly, we obtain a double stability theorem for Serrin's overdetermined problem in spaceforms. Secondly, we prove stability theorems for Brendle's Heintze-Karcher inequality respectively constant mean curvature classification in a class of warped product spaces. The key tool is the first author's recent development of stability for level sets of a function under smallness of the traceless Hessian thereof.

\end{abstract}

\maketitle

\section{Introduction}

This paper is about three famous problems from geometric analysis, which are deeply interrelated. The study of constant mean curvature (CMC) hypersurfaces has a very long history. While being far from a full classification in general, for example in the complete or immersed class, the case of closed and embedded CMC hypersurfaces of the Euclidean space is fully understood, due to Alexandroff \cite{Alexandroff:12/1962}: All of them are round spheres. This result also holds true in the spherical and hyperbolic spaces. In more general warped spaces, such as the Schwarzschild spaces, results of this type are more recent and due to Brendle \cite{Brendle:06/2013}. While Alexandroff's original proof is an ingenious application of the elliptic maximum principle and a reflection technique, rather different proofs have emerged later on. For example, inspired by an idea of Reilly \cite{Reilly:/1977b}, Ros \cite{Ros:/1987} gave a very elegant proof via integral methods. He proved an inequality which by today is known as Heintze-Karcher inequality (HK), and which reads
\eq{\int_{\del\Om}\fr{n}{H}\geq (n+1)\abs{\Om} = \int_{\del\Om}\ip{x}{\nu},}
where $\Om\sub\bbR^{n+1}$ is a domain with volume $\abs{\Om}$ and mean-convex boundary, $\nu$ is the outward pointing normal, $x$ the position and $H$ the trace of the second fundamental form of the boundary a.k.a mean curvature. The key in this inequality is that equality occurs if and only if $\del\Om$ is a round sphere. So if $H$ is constant, then 
\eq{\int_{\del\Om}\ip{x}{\nu} = \frac{1}{H} \int_{\del\Om}H\ip{x}{\nu} = \fr{n\abs{\del\Om}}{H} = \int_{\del\Om}\frac{n}{H},}
where in the middle identity we used
\eq{0=\int_{\del\Om}\dive_{M}x^{\top} = \int_{\del\Om}(n-H\ip{x}{\nu})}
and where $\abs{\del\Om}$ is the surface area of $\del\Om$.
Hence equality in HK holds, and thus $\del\Om$ must be a round sphere. CX showed together with Qiu \cite{QiuXia:/2015} and respectively with Li \cite{LiXia:11/2019} that there also is such a type of proof for Brendle's HK respectively CMC classification in spaceforms and more general warped products. Also compare the refinement \cite{FogagnoloPinamonti:07/2022}. Our third problem of interest is Serrin's overdetermined problem. Note that one way to prove HK is to solve the so-called {\it{torsional rigidity equation}}
\eq{\label{torsion}\bar\De f &= 1\q\mbox{in}~\Om\\
		f&=0\q\mbox{on}~\del\Om,}
where $\bar\De$ is the standard Laplace operator of the ambient Euclidean space. In case $\Om$ is a ball, the function $f$ is, up to normalisation and centre, the squared distance function to a point, in which case the outward Neumann derivative $f_{\nu}$ of $f$ is constant. Indeed, the reverse is also true and this result is the famous overdetermined problem of Serrin \cite{Serrin:01/1971}: If we add the condition $f_{\nu} = \const$ to \eqref{torsion}, then $\Om$ is a ball and $f$ is a squared distance function. A vast amount of {\it{Serrin-type}} results have been achieved in various settings, for example in spaceforms, \cite{CiraoloVezzoni:07/2019,KumaresanPrajapat:/1998,QiuXia:/2017}.

The characterisation of the domains satisfying the respective equality cases as being geodesic balls in the above problems, namely 
\eq{H=\const,\q \int_{\del\Om}\fr{n}{H} = (n+1)\abs{\Om},\q f_{\nu} = \abs{\bar\n f} = \const,}
allows to ask the associated {\it{stability question}}: Must $\del\Om$ be close to a sphere, provided that one of the above equalities are ``almost'' satisfied? For the CMC problem this question was affirmatively answered in great detail in spaceforms and with various different methods, \cite{CiraoloVezzoni:02/2018,CiraoloVezzoni:/2020,MagnaniniPoggesi:10/2019}, while the stability for HK was treated in \cite{MagnaniniPoggesi:10/2019,Scheuer:03/2021}. The stability for Serrin's problem was considered in \cite{MagnaniniPoggesi:01/2020} and we refer to \cite{Scheuer:03/2021} for a deeper discussion. In the latter paper, JS achieved what was called {\it{double stability result}} for Serrin's problem in the Euclidean space, which refers to the fact that not only the normal derivative is allowed to be perturbed, but also the right hand side of the Poisson equation. One of the three results we want to present here is a generalisation of this double stability result to spaceforms. In particular we are going to prove the following theorem, where $\bbM_{K}$ is the hyperbolic space for $K=-1$ or an open half-sphere in case $K=1$ with respective origin $O$, and where $\abs{\cdot}_{k,\al,\Om}$ and $\|\cdot\|_{p,\Om}$ denote the $C^{k,\alpha}$- resp. the $L^{p}$-norm on $\Om$,
\eq{ \|f\|_{p,\Om} = \br{\int_{\Om}\abs{f}^{p}}^{\fr 1p}, \quad \|f\|_{\infty, \Om}=\sup_{\Om} |f|, }
and similarly for tensors on manifolds. $\bar\De$ is the Laplace operator of the metric on the ambient space and 
\eq{V = \vt'(r),}
where $r$ is the distance function to the origin of the spaceform and 
\eq{\vt(r) = \begin{cases} \sin r, & K=1\\
					\sinh r, & K=-1.\end{cases}} 
Furthermore, $A$ is the second fundamental form on $\del\Om$ and $\mr{A}$ its traceless part.

\begin{thm}[Double stability in Serrin's problem]\label{main thm}
Let $n\geq 2$, $K=\pm 1$, $\Om\sub \bbM_{K}^{n+1}$ be a bounded domain with connected $C^{2,\beta}$-boundary ($0<\be < 1$) and uniform interior ball condition with radius $\rho$ and $\bar\Om\sub \bbM_{K}$.
Let $\phi\in C^{1}((-\8,0])$ be positive and suppose that $f\in C^{2,\be}(\bar\Om)$ satisfies
\eq{\label{equation}\bar\De f +(n+1)Kf &= \phi(f),\q&&\mbox{in}~\Om\\
				f&=0\q&&\mbox{on}~\del\Om.}
Then there exists
\eq{C = C(n,\abs{f}_{2,\be,\Om},\|\phi\|_{\8,\Om},\|1/\phi\|_{\8,\Om},\rho^{-1},\abs{\del\Om},\max_{\bar \Om}r),\footnotemark}
\footnotetext{For $K=1$, the dependence on $\max r$ is understood to become singular at $\pi/2$.}such that
\eq{\|\mr{A}\|^2_{2,\del\Om}&\leq C(\|f_{\nu}-R\|_{1,\del\Om}+\|\phi-1\|_{1,\Om} + \|\Phi-f\phi\|_{1,\Om}+\|f\Phi-\Psi-\tfr 12 f^{2}\|_{1,\Om})^\frac{\be}{2+\be},
}
where
\eq{R=\fr{\int_{\Om}V\phi}{\int_{\del\Om}V} ,\q\Phi(f) = \int_{0}^{f}\phi(s)\,ds\q \mbox{and}\q \Psi(f) = \int_{0}^{f}\Phi(s)\,ds.}
\end{thm}

\begin{rem}
From \autoref{main thm}, it is possible to deduce various kinds of closeness of $\del\Om$ to a geodesic sphere, depending on what additional regularity one is given for $f$. For example, with the given regularity one can deduce for $p>n$,
\eq{\|\mr{A}\|_{p,\del\Om}^{p}\leq C\|\mr{A}\|^{2}_{2,\del\Om}}  
and hence it is possible to obtain $W^{2,p}$-closeness of $\del\Om$ to a geodesic sphere by using \cite[Thm.~1.3]{De-RosaGioffre:/2021}, after using a suitable conformally flat parametrisation.  But there are many other stability results, such as the results of De Lellis/M\"uller \cite{De-LellisMuller:/2005,De-LellisMuller:/2006}.
\end{rem}

As mentioned above, stability results for the CMC problem and for HK in spaceforms are known. Here we would like to extend these results to a class of warped product spaces, where in particular it appears to be of interest to obtain quantitative versions of the corresponding rigidity results of Brendle \cite{Brendle:06/2013}. For this purpose, we prove the following two results, where $H_{1} = H/n$ is the normalised mean curvature, and where $\abs{M}_{2,\beta_{2}}$ is an abbreviation for control of the geometry of $M$ in $C^{2,\beta_{2}}$-sense, as it appears in the constant of the Schauder estimates, see \eqref{reg-f}.

\begin{thm}[Stability of HK in warped products]\label{main:HK}
Let 
\eq{\bar M^{n+1} = [0,\overline{r})\x\bbS^{n},\q\bar g = dr^{2}+\vt^{2}(r)\si} be a warped product manifold with $\vt$ satisfying \eqref{H1}-\eqref{H3} and \eqref{H5}, where $\si$ is the round metric. Let $0<\be_{2}<1$ and $M\sub \bar M$ be a strictly mean convex $C^{2,\beta_2}$-hypersurface, which, possibly in conjunction with $\{r=0\}$, bounds a domain $\Om$
 with uniform interior ball condition of radius $\rho$. Then 
  \eq{C=C\br{n, \rho^{-1}, |M|, \max_M r,  (\min_M r)^{-1},|M|_{2,\beta_2,}},}
 such that
\eq{
\|\mr{A}\|^2_{2, M}\le  C\br{\int_M \left(\frac{\vt'}{H_{1}} -\bar g(\vt\del_r, \nu)\right)}^{\frac{\beta}{1+\beta}}.
}
where $\beta=\min\{\beta_1,\beta_2\}$.
\end{thm}

For the additional notation involved here, please consult \autoref{Prelim}. The properties \eqref{H1}-\eqref{H5} resemble closely the properties in \cite{Brendle:06/2013} and the details are given at the beginning of \autoref{CMC}.
We remark that, since $\{r=0\}\times \bbS^n$ is minimal and $M$ is strictly mean convex, $M$ cannot touch $\{r=0\}\times \bbS^n$ by the maximum principle. Hence $\min_M r>0$. In the case that $M$ is homologous to $\{r=0\}\times \bbS^n$, $(\min_M r)^{-1}$ can be controlled by $\rho^{-1}$.

Finally, very much related to this result is the stability of the CMC problem, meaning that every almost-CMC hypersurface in warped products satisfying \eqref{H1}-\eqref{H5} must be close to a coordinate slice.

\begin{thm}\label{main:CMC}
Let $(\bar M, \bar g)$ be a warped product manifold with $\vt$ as in \autoref{main:HK} satisfying \eqref{H1}-\eqref{H5}. Let $M\sub (0,\bar r)\x\bbS^{n}$ be a  $C^{2,\beta_2}$-hypersurface, which, possibly in conjunction with $\{r=0\}$, bounds a domain $\Om$ with uniform interior ball condition of radius $\rho$. 
Then there exist constants $\al(n)$
 and $C,\be$ as in  \autoref{main:HK},
such that
\eq{\label{cmc-stab}
d(M, \{r=r_0\}\times \mathbb{S}^n)\le C\left\|H_1-\frac{\int_{M}\vt'}{(n+1)\int_{\Om}\vt'+\vt(0)^{1+n}|\bbS^n|}\right\|_{1,M}^{\alpha(n)\beta}
}
for some $r_0\in (0, \bar r)$. Here $d$ is the Hausdorff distance.
\end{thm}

Recall that Li and the second author \cite{LiXia:11/2019} reproved Brendle's HK and CMC classification by applying the solution to an appropriate boundary value problem in a new Reilly type formula. The main idea for the stability results in \autoref{main:HK} and \autoref{main:CMC} is to provide this proof in a quantitative way, with the help of a key idea due to the first author \cite{Scheuer:03/2021}.

\begin{rem}
From the proof of \autoref{main:CMC}, in particular \eqref{pf:thm-CMC-1}, it can be seen that in case the right hand side of \eqref{cmc-stab} is small, the hypersurface $M$ is actually graphical over $\bbS^{n}$ and we have an estimate of the $W^{1,p}$-norm in terms of that quantity. Higher regularity can then be deduced, depending on the regularity of the curvature. We do not address these questions further, as they can be derived in standard ways from interpolation inequalities.
\end{rem}

\section{Preliminaries}\label{Prelim}
Let $\bbS^{n}$ be the $n$-dimensional unit sphere equipped with the round metric $\si$. For a real number $\overline{r}>0$, let $\bar{M}^{n+1}=[0,\overline r)\times \bbS^n$ be a warped product manifold with metric 
\eq{\bar g=dr^2+\vt^{2}(r) \sigma,}
where $\vt: [0,\overline r)\to \bbR_+$ is a smooth function.

In the simply connected spaceforms with sectional curvature $K$ there holds 
\eq{\vt(r) = \begin{cases} \sin r,& K=1\\
						r,& K=0\\
					\sinh r,& K=-1.\end{cases}}
We denote by $\bar\n$, $\bar\De$ and $\ov{\Rc}$ the Levi-Civita connection, the Laplace operator and the Ricci tensor of $\bar g$, respectively.
By direct computation as in \cite{Brendle:06/2013,BrendleEichmair:07/2013}, we see that
\eq{\label{Rc} \ov{\Rc}&=\left((n-1)\frac{1-\vt'^2}{\vt^2}-\frac{\vt''}{\vt}\right)\bar g-(n-1)\left(\frac{\vt''}{\vt}+\frac{1-\vt'^2}{\vt^2}\right)dr\otimes dr 
\\&=\left((n-1)(1-\vt'^2)-\vt\vt''\right)\sigma-n\frac{\vt''}{\vt}dr\otimes dr.
}
When $\vt'>0$, we also have
\eq{\label{D^{2}theta'}
&\frac{\bar \n^2 \vt'}{\vt'}=\vt\vt''\sigma+\frac{\vt'''}{\vt'}dr\otimes dr}
and
\eq{\frac{\bar \De \vt'}{\vt'}=\frac{\vt'''}{\vt'}+n\frac{\vt''}{\vt}.
}
Thus 
\eq{\label{Laplace theta'}
&\frac{\bar \De \vt'}{\vt'}\bar g- \frac{\bar \n^2 \vt'}{\vt'}+\ov{\Rc} =\left(\vt^2\frac{\vt'''}{\vt'}+(n-2)\vt\vt''+(n-1)(1-\vt'^2)\right)\sigma.
}

It is known that the vector field $X=\vt(r)\del_{r}$ is a conformal Killing vector field such that the Lie derivative satisfies  
 \eq{L_X\bar g=2\vt'(r)\bar g,} and hence the function
\eq{V = \vt'(r)}
plays a very special role in the subsequent considerations. 
For a two-tensor $T$ denote by
\eq{]T[_{\al\be} = \tfr 12(T_{\al\be}-T_{\be\al})} its anti-symmetri\-sation
and by
\eq{[T]_{\al\be} = \tfr 12(T_{\al\be}+T_{\be\al})}
its symmetrisation.
The conformal Killing property of $X$ can then also be expressed as
\eq{[\bar\n X^{\flat}] = V\bar g,}
where $X^{\flat}$ is defined in the standard way by 
\eq{X^{\flat}(Y) = \bar g(X,Y)}
for all tangent vectors $Y$.
Furthermore, we obtain for any $C^{2}$-function $f$,
\eq{\bar\n\br{V^{2}\bar\n\br{\fr fV}} = \bar\n\br{V\bar\n f - f\bar\n V} = V\bar\n^{2}f - f\bar\n^{2}V + 2]\bar\n V\otimes\bar\n f[}
and hence
\eq{\label{divergence}\left[\bar\n\br{V^{2}\bar\n\br{\fr fV}}\right] = Vf\br{\fr{\bar\n^{2}f}f-\fr{\bar\n^{2}V}{V}}.}
It follows by taking trace of \eqref{divergence} that
\eq{\label{divergence1}\ov{\dive}\br{V^{2}\bar\n\br{\fr fV}}= Vf\br{\fr{\bar\De f}{f} - \fr{\bar\De V}{V}}.}
The divergence structure in \eqref{divergence} and \eqref{divergence1} has been noticed in \cite{LiXia:11/2019}, see also \cite{LiXia:07/2017} for a more general setting.

In spaceforms of sectional curvature $K$, these formulas simplify a bit because in this case
\eq{\label{Hessian V}\bar\n^{2}V = -KV\bar g.}
Using the symbol $\bar g$ to denote the natural inner product induced on all higher tangent- and cotangent bundles, we can write
\eq{\label{Codazzi}V\abs{\bar\n^{2}f+Kf\bar g}^{2} &= \bar g\br{\bar\n\br{V^{2}\bar\n\br{\fr{f}{V}}},\bar\n^{2}f+Kf\bar g}\\
					&=\ov{\dive}\br{(\bar\n^{2}f + Kf\bar g)\br{V^{2}\bar\n\br{\fr{f}{V}}}}\\
					&\hp{=} - V^{2}\bar g\br{\bar\n\br{\fr{f}{V}},\bar\n(\bar\De f+(n+1)Kf)},}
where we have used that $\bar\n^{2}f+Kf\bar g$ is Codazzi.

From now on, we simplify notation by defining
\eq{\ip{Y}{Z} := \bar g(Y,Z)}
for vector fields $Y,Z$. We write $\bbM_{K}$ for the simply connected spaceform with sectional curvature $K$, where in the spherical case $K=1$, we put $\overline{r} = \pi/2$. 
For a bounded domain $\Om\sub\bar M$ we use $\nu$ as our default symbol for the outward pointing unit normal. This means that at a maximal point of the radial coordinate $r$ we always assume $\nu = \del_{r}$. For a function $f\in C^{1}(\bar\Om)$, the Neumann derivative in direction $\nu$ is written $f_{\nu}$. If $M\sub \bar M$ is a $C^{2}$-hypersurface, we define the second fundamental form $A=(h_{ij})$ by
\eq{\bar\n_{Y}Z = \n_{Y}Z - A(Y,Z)\nu, }
where $\n$ is the Levi-Civita connection of the induced metric $g$ in $M$.
For a two-tensor $T$ on any Riemannian manifold, $\mr{T}$ denotes its tracefree part.
 In $\bbM_{K}$, we denote by $O$ the origin, where $r=0$. For any $k$-dimensional closed submanifold $\Si$, $\abs{\Si}$ denotes the volume associated with the induced volume element.  Finally, constants $C = C(\dots)$ may vary from line to line, as long as no other than the indicated dependencies are introduced.


\section{Stability for Serrin's problem with nonlinear source in spaceforms	}\label{Serrin}

\subsection*{A Pohozaev type identity}

\begin{lemma}\label{Pohozaev}
Let $\Om\sub\bbM_{K}$ be a $C^{1}$-domain and let $f\in C^{2}(\bar\Om)$ be a solution to \eqref{equation}, then
\eq{\fr{1}{2}\int_{\del\Om}\ip{X}{\nu}f_{\nu}^{2} &= \fr{(n+1)(n+3)K}{4}\int_{\Om}Vf^{2} - (n+1)\int_{\Om} V\Phi \\
						&\hp{=} +\fr{n-1}{2}\int_{\Om}Vf\phi(f),}
where $X=\vt(r)\del_{r}$.
\end{lemma}

\pf{
We calculate on the one hand,
\eq{\int_{\Om}\ip{X}{\bar\n f}\bar\De f &= \int_{\del\Om}\ip{X}{\bar\n f}f_{\nu} - \int_{\Om}d\ip{X}{\bar\n f}(\bar\n f)\\
						&= \int_{\del\Om}\ip{X}{\bar\n f}f_{\nu} - \int_{\Om}V\abs{\bar\n f}^{2} - \fr 12\int_{\Om}d(\abs{\bar\n f}^{2})(X)\\
						&= \fr 12\int_{\del\Om}\ip{X}{\nu}f_{\nu}^{2} + \fr{n-1}{2}\int_{\Om}V\abs{\bar\n f}^{2},}
and on the other hand,
\eq{\int_{\Om}\ip{X}{\bar\n f}\bar\De f &= \int_{\Om}\ip{X}{\bar\n f}(\phi(f) - (n+1)Kf)\\
							&=\int_{\Om}\ip{X}{\bar\n\br{\Phi(f) - \fr{(n+1)K}{2}f^{2}}}\\
							&=\int_{\Om}\fr{(n+1)^{2}KV}{2}f^{2} - \int_{\Om}(n+1)V\Phi.}
Finally,
\eq{\int_{\Om}V\abs{\bar\n f}^{2}&=\int_{\Om}V\br{\tfr 12\bar\De(f^{2})-f\bar\De f}\\
							&=\fr 12\int_{\Om} f^{2}\bar\De V +  \int_{\del\Om}Vf\ip{\bar\n f}{\nu} - \fr 12\int_{\del\Om}f^{2}\ip{\bar\n V}{\nu}\\
								&\hp{=}-\int_{\Om}fV(\phi(f)-(n+1)Kf)\\
							&=\fr{(n+1)K}{2}\int_{\Om} f^{2}V -\int_{\Om}fV\phi(f).}
The proof is complete by combining these three equalities.
}

\subsection*{Proof of \autoref{main thm}}
 We aim to make use of the key idea from \cite{Scheuer:03/2021} and therefore we are going to provide an estimate for the quantity
\eq{\int_{\Om}Vf\abs{\mr{\bar\n}^{2}f}^{2} &= \int_{\Om}Vf\abs{\bar\n^{2}f - \tfr{1}{n+1}\bar\De f \bar g}^{2}\\
	& =\int_{\Om}Vf\abs{\bar\n^{2}f+Kf\bar g - \tfr{1}{n+1}\br{\bar\De f + (n+1)Kf}\bar g}^{2}.}
Using \eqref{equation}, \eqref{divergence}, \eqref{Codazzi}  and multiple integration by parts, we obtain
\eq{&\int_{\Om}Vf\abs{\mr{\bar\n}^{2}f}^{2}\\
	 = &\int_{\Om}Vf\abs{\bar\n^{2}f + Kf\bar g}^{2}-\fr{1}{n+1}\int_{\Om}Vf\phi^{2} \\
							= &\int_{\Om}f\ov{\dive}\br{(\bar\n^{2}f + Kf\bar g)\br{V^{2}\bar\n\br{\tfr{f}{V}}}} - \int_{\Om}fV^{2}\ip{\bar\n\br{\tfr{f}{V}}}{\bar\n\phi}-\fr{1}{n+1}\int_{\Om}Vf\phi^{2}\\
					= &-\int_{\Om} (\bar\n^{2}f+Kf\bar g)\br{\bar\n f,V^{2}\bar\n\br{\tfr{f}{V}}}-\int_{\Om}\ip{V^{2}\bar\n\br{\tfr{f}{V}}}{\bar\n(f\phi-\Phi)} \\
					&- \fr{1}{n+1}\int_{\Om}Vf\phi^{2}\\
					=&-\int_{\Om}\tfr{1}{2}\ip{V^{2}\bar\n\br{\tfr{f}{V}}}{\bar\n\abs{\bar\n f}^{2}+K\bar\n f^{2}}+\int_{\Om}V\phi(f\phi-\Phi)  -\fr{1}{n+1}\int_{\Om}Vf\phi^{2}\\
					= &\int_{\Om}V\phi(\tfr 12 \abs{\bar\n f}^{2}+ \tfr 12 Kf^{2}+\tfr{n}{n+1}f\phi - \Phi)  - \fr 12\int_{\del\Om}Vf_{\nu}^{3}. }

Now we use
\eq{\fr 12\int_{\Om}V\phi\abs{\bar\n f}^{2} &= -\fr 12\int_{\Om}\Phi\ip{\bar\n V}{\bar\n f} - \fr 12\int_{\Om}V\Phi\bar\De f\\
					&= -\fr 12\int_{\Om}\ip{\bar\n V}{\bar\n \Psi} - \fr 12\int_{\Om}V\Phi(\phi - (n+1)Kf)\\
					&=  \fr{(n+1)K}{2}\int_{\Om}V(\Phi f-\Psi) - \fr 12\int_{\Om}V\Phi\phi}
and plug in \autoref{Pohozaev} in order to compare the boundary term:

\eq{&\int_{\Om}Vf\abs{\mr{\bar\n}^{2}f}^{2}\\
	=&\fr{(n+1)K}{2}\int_{\Om}V(\Phi f-\Psi) - \fr 12\int_{\Om}V\Phi\phi+\int_{\Om}V\phi( \tfr 12 Kf^{2}+\tfr{n}{n+1}f\phi - \Phi)\\
	&  - \fr 12\int_{\del\Om}(Vf_{\nu} - \tfr{1}{n+1} \ip{X}{\nu})f_{\nu}^{2}-\fr{1}{2(n+1)}\int_{\del\Om}\ip{X}{\nu}f_{\nu}^{2}\\
	=&\fr{(n+1)K}{2}\int_{\Om}V(\Phi f-\Psi) - \fr 32\int_{\Om}V\Phi\phi+\int_{\Om}V\phi( \tfr 12 Kf^{2}+\tfr{n}{n+1}f\phi)\\
	&  - \fr 12\int_{\del\Om}(Vf_{\nu} - \tfr{1}{n+1} \ip{X}{\nu})(f_{\nu}^{2}-R^{2}) - \fr{R^{2}}{2}\int_{\del\Om}(Vf_{\nu} - \tfr{1}{n+1}\ip{X}{\nu})\\
	&-\fr{(n+3)K}{4}\int_{\Om}Vf^{2} + \int_{\Om} V\Phi -\fr{n-1}{2(n+1)}\int_{\Om}Vf\phi.}

Now use
\eq{\int_{\del\Om}Vf_{\nu} = \int_{\Om}V\phi\q\mbox{and}\q\int_{\del\Om}\ip{X}{\nu} = -K\int_{\del\Om}V_{\nu} = (n+1)\int_{\Om}V}
to obtain
\eq{\label{pf:Serrin-1}&\int_{\Om}Vf\abs{\mr{\bar\n}^{2}f}^{2}\\
	=&\fr{R^{2}}{2}\int_{\Om}V(1-\phi)  - \fr 12\int_{\del\Om}(Vf_{\nu} - \tfr{1}{n+1} \ip{X}{\nu})(f_{\nu}^{2}-R^{2}) \\
	&+\fr{(n+1)K}{2}\int_{\Om}V(\Phi f-\Psi) -\fr{(n+3)K}{4}\int_{\Om}Vf^{2}+\fr{K}{2}\int_{\Om}V\phi f^{2}\\
	& + \fr{n}{n+1}\int_{\Om}Vf\phi^{2}- \fr 32\int_{\Om}V\Phi\phi + \int_{\Om} V\Phi -\fr{n-1}{2(n+1)}\int_{\Om}Vf\phi\\
	=&\fr{R^{2}}{2}\int_{\Om}V(1-\phi)  - \fr 12\int_{\del\Om}(Vf_{\nu} - \tfr{1}{n+1} \ip{X}{\nu})(f_{\nu}^{2}-R^{2}) \\
	&+\fr{(n+1)K}{2}\int_{\Om}V(\Phi f-\Psi -\tfr 12 f^{2}) +\fr{K}{2}\int_{\Om}V(\phi-1) f^{2}\\
	& -\fr{n+3}{2(n+1)}\int_{\Om}Vf\phi(\phi-1) + \int_{\Om}V(1-\tfr 32\phi)(\Phi-f\phi) .}

From here, the proof can be completed in a similar way as the proof of \cite[Thm.~1.10]{Scheuer:03/2021}, after we also invoke the additional regularity of $f$.
From \cite[Lemma~4.1]{Scheuer:03/2021} we obtain $f<0$ in $\Om$ and control on $\abs{\bar\n f}$ from below in terms of $\min\phi$ and $\rho$, the radius satisfying the interior sphere condition for $\Om$. We define
\eq{\ep = \|f_{\nu}-R\|_{1,\del\Om}+\|\phi-1\|_{1,\Om} + \|\Phi-f\phi\|_{1,\Om}+\|f\Phi-\Psi-\tfr 12 f^{2}\|_{1,\Om}}
and may wlog assume that $\ep>0$, for otherwise we are in the classical case of Serrin's theorem. We also assume $\ep<1$, for otherwise the claimed estimate is trivial due to the allowed dependencies in $C$.

Using the flow $F$ defined by 
\eq{\label{level-set-flow}\del_{t}F(t,\xi) &= -\fr{\bar\n f(F(t,\xi))}{\abs{\bar\n f(F(t,\xi))}^{2}}\\
				F(0,\xi)& = \xi,}
 we can construct a one-sided open neighbourhood of $M$ of the form
\eq{\label{level-set-nbhd}\cU = \bigcup_{0<t<T}\{f_{|\cU}=-t\}\equiv \bigcup_{0<t<T}M_{-t},}
where  
\eq{T = \min\Big\{\fr{\min_{M}\abs{\bar\n f}^{2}}{4\|\bar\n^{2}f\|_{\8,\Om}},\ep^{\fr{1}{2+\be}}\Big\}.}
From the mean value theorem we obtain on $\cU$,
\eq{\abs{\bar\n f}\geq \tfr12\min_{M}\abs{\bar\n f}}
and furthermore the level surfaces $M_{-t}$ are connected. Fixing $t$, then for a local frame $(e_{i})$ of $M_{-t}$, from the relation
\eq{\bar\n^{2}f(e_{i},e_{j}) = -\abs{\bar\n f}h_{ij}}
a simple calculation gives
\eq{\abs{\bar\n f}^{2}\abs{\mr{A}}^{2}=\mr{\bar\n}^{2}f(e_{i},e_{j})\mr{\bar\n}^{2}f(e_{k},e_{l})g^{ik}g^{jl}-\tfr 1n (\mr{\bar\n}^{2}f(\nu,\nu))^{2},}
where $\nu=\bar\n f/\abs{\bar\n f}$ on $M_{t}$. The co-area formula and \eqref{pf:Serrin-1} give
\eq{\label{co-area}\int_{-T}^{-T/2}\int_{M_{s}}\abs{\bar\n f}^{}\abs{\mr{A}}^{2}\,ds\leq \int_{\{-T<t<-T/2\}}\abs{\mr{\bar\n}^{2}f}^{2}\leq \fr{C}{T}\int_{\Om}V(-f)\abs{\mr{\bar\n}^{2}f}^{2}\leq \fr{C\ep}{T}.}
We can ignore the gradient term on the left hand side, because it is under control from below. Hence there exists $s\in[-T,-T/2]$, such that
\eq{\|\mr{A}\|^{2}_{2,M_{s}}\leq \fr{C\ep}{T^{2}}.}
Finally, we have to transfer this estimate to $M_{0}=\del\Om$ and therefore we use the H\"older regularity of $\bar\n^{2}f$. Fix $\xi\in M_{0}$ and let $(t,\xi)$ denote the coordinates of $\cU$, then
\eq{\abs{\abs{\mr{A}(s,\xi)}^{2}-\abs{\mr{A}(0,\xi)}^{2}}\leq C(-s)^{\be}\leq CT^{\be},}
where we used that the flow $F$ is of class $C^{1,\be}$, due to the regularity of $f$, cf. \cite[Thm.~9.4.4]{Gerhardt:/2006b}. Hence
\eq{\left|\int_{M_{s}}\abs{\mr{A}}^{2}-\int_{M_{0}}\abs{\mr{A}}^{2}\right|\leq CT^{\be}}
and
\eq{\|\mr{A}\|^{2}_{2,\del\Om}\leq \fr{C\ep}{T^{2}}+CT^{\be}\leq C\ep^{\frac{\be}{2+\be}}.}
The proof is complete.

\section{Stability for the CMC problem and the Heintze-Karcher inequality in a class of warped products}\label{CMC}

 We assume the warped product factor $\vt\cn [0,\bar r)\to \bbR_+$ is a smooth positive function which satisfies the following:
\begin{align}\label{H1}\tag{H1} &\vt'(0)=0,\quad \vt''(0)>0,\\
\label{H2}\tag{H2}&\vt'(r)>0 \hbox{ for }r\in (0,\bar r), \\
\label{H3}\tag{H3}& 2\frac{\vt''}{\vt}-(n-1)\frac{1-\vt'^2}{\vt^2}\hbox{ is non-decreasing for }r\in (0,\bar r),\\
\label{H4}\tag{H4}& \frac{\vt''}{\vt}+\frac{1-\vt'^2}{\vt^2}>0, \hbox{ for }r\in (0,\bar r),\\
\label{H5}\tag{H5}& \frac{\vt'''}{\vt'}\in C^{0,\be_{1}}([0,\bar r)) \hbox{ for some }\beta_1\in (0, 1).
\end{align}

Assumptions \eqref{H1}-\eqref{H4} coincide with those in Brendle \cite{Brendle:06/2013}, while we need an additional assumption \eqref{H5}, because we need to solve an elliptic boundary value problem \eqref{dirichlet} below in which $\bar \De \vt'/\vt'$ occurs as the coefficient for an elliptic operator. Since $\vt'(0)=0$, \eqref{H5} says that $\vt'''/\vt'$ is H\"older continuous up to $\{r=0\}$.

Typical examples which satisfy \eqref{H1}-\eqref{H5} include the (de-Sitter-)Schwarzschild manifolds and the Reissner-Nordstrom manifolds.  Conditions \eqref{H1}-\eqref{H4} for such examples\footnote{(H4) holds when $m>0$.} were verified in Brendle \cite{Brendle:06/2013}. Here we verify \eqref{H5} for such examples.
The warped factor $\vt$ of (de-Sitter-)Schwarzschild manifold solves 
\eq{
\vt'(r)=\sqrt{1+\kappa \vt^2-2m \vt^{1-n}}, 
}for some $\kappa=1, 0\hbox{ or } -1$ and $m\in \bbR$.
Thus one computes that 
\eq{
\frac{\vt'''}{\vt'}=\kappa - {mn(n-1)\vt^{-n-1} }, 
}
which is readily smooth in $[0, \bar r)$. 
Similarly, for the Reissner-Nordstrom manifold, where \eq{
\vt'=\sqrt{1+\kappa \vt^2-2m \vt^{1-n}+q^2\vt^{2-2n}}, 
}for some $\kappa=1, 0\hbox{ or } -1$,  $m\in \bbR$ and $q>0$, one gets that $\vt'''/\vt'$ is also smooth in $[0, \bar r)$.

Note that by \eqref{H1} the slice $\{r=0\}\times \bbS^n$ is totally geodesic, since 
\eq{ h_{ij}|_{\{r=0\}\times \bbS^n}=\bar g(\bar \n_{\del_i}\del_r, \del_j)= \vt'(0)\vt(0)\sigma_{ij}=0.
}

Equation \eqref{Laplace theta'} implies that \eqref{H3} is equivalent to 
\eq{\bar \De \vt'\bar g-\bar \n^2 \vt'+\vt'\ov{\Rc}\ge 0.}
\eqref{H5} implies that $\bar \De \vt'/\vt'$ is $\beta_1$-H\"older continuous up to $\{r=0\}\times \bbS^n$.

\subsection*{Stability in the Heintze-Karcher inequality}

We need the following Reilly type formula due to Li and the second author \cite{LiXia:11/2019}.
\begin{prop}[\cite{LiXia:11/2019}, Theorem 1.1]
Let $\Om$ be a bounded domain in $(\bar M, \bar g)$ with boundary $\del \Om$. 
 Let $f, V\in C^2(\bar \Om)$ such that $\frac{\bar \nabla^2 V}{V}$ is continuous up to $\del \Om$. Then there holds
\eq{\label{lixia}
           &\int_{\Om}V \Big(\bar \Delta f-\frac{\bar \Delta V}{V}f\Big)^2-\int_{\Om}V\Big|\bar \nabla^2 f-\frac{\bar \nabla^2 V}{V}f\Big|^2 
           \\=&\int_{\del\Om}V h(\n f,\n f)+2V f_\nu\Delta f+nV H_1f_\nu^2+V_{\nu}|\n f|^2+2f\bar \nabla^2 V(\n f, \nu) dA\\ 
            &+\int_{\del\Om} -2ff_\nu (\De V+nH_1V_{\nu})-f^2\frac{\bar \nabla^2 V-\bar \Delta V \bar g}{V}(\bar \nabla V, \nu) dA\\
            &+  \int_{\Om}(\bar \Delta V g- \bar \nabla^2 V+V\ov{\Rc})\left(\bar \nabla f-\frac{\bar \nabla V}{V}f, \bar \nabla f-\frac{\bar \nabla V}{V}f\right).
}
\end{prop}

\subsection*{Proof of \autoref{main:HK}}
Let $V=\vt'(r)$. Let $\Om$ be the bounded domain in $\bar M$ enclosed by $M$ in the case  $M$ is null-homologous and by $M$ and $\{r=0\}\times \bbS^n$ in the case  $M$ is homologous to $\{r=0\}\times \bbS^n$. Consider the Dirichlet problem 
\eq{\label{dirichlet} &\bar \De f-\frac{\bar \De V}{V} f=1&&\hbox{ in }\Om,\\
&f=0&&\hbox{ on }M,\\
&f=c_0&&\hbox{ on }\{r=0\}\times \bbS^n,
}
where the final boundary condition shall be ignored in case that $M$ is null-homo\-logous and where \begin{eqnarray}\label{c0}
c_0= - \frac{n}{n+1}\frac{V}{(\bar \Delta V \bar g-\bar \nabla^2 V)(\nu, \nu)}\Big|_{\{r=0\}\times \bbS^n}=-\frac{1}{n+1} \frac{\vt(0)}{\vt''(0)}.
\end{eqnarray}

Using Moser's iteration to estimate the $L^\infty$ norm of $f$  by its $L^{2}$ norm, see \cite[Thm.~8.15]{GilbargTrudinger:/2001}, a corollary of the Fredholm alternative, see \cite[Cor.~8.7]{GilbargTrudinger:/2001}, and Schauder theory, we obtain the existence of a unique solution $f\in C^{2,\beta}(\bar \Om)$ to \eqref{dirichlet}, where $\beta=\min\{\beta_1,\beta_2\}$ with the estimate  
\eq{\label{reg-f}|f|_{2,\be,\Om}\le C\left(|M|_{2,\beta_2}, \left|\frac{\bar \De V}{V}\right|_{0,\beta_1, \Om}\right).
}
The crucial point to use \cite[Cor.~8.7]{GilbargTrudinger:/2001} is that zero is not an eigenvalue of $L=\bar \Delta-\frac{\bar \Delta V}{V}$, as shown in \cite[Lemma~2.5]{LiXia:11/2019}.

Next we prove that $f<0$ in $\Om$. As in \eqref{divergence1}, 
we can transform the equation \eqref{dirichlet} into
\eq{\label{dirichlet-2}\ov{\dive}\left(V^2\bar \n\br{\fr fV}\right)=V.
}

It follows from the strong maximum principle that $\frac{f}{V}$ cannot attain an interior maximum. From \eqref{H1}, we see $c_0<0$. Hence $\frac{f}{V}<0$ and in turn $f<0$ in $\Om$.
By using the standard barrier argument, one deduces a  positive lower bound for the boundary gradient:
\eq{\min_{M}\abs{\bar\n f}\ge C>0,}
where $C$ depends on  $\rho^{-1}$, $\max_{M} r$ and $\min_{M}r$.
The reason for this can be given as follows: Pick a boundary point $x\in M$ and let $B = B_{\rho/2}(x_{0})$ be an interior ball. The warped product is diffeomorphic to a Euclidean annulus and hence we can write \eqref{dirichlet-2} within $B$ in $\bar g$-conformal Euclidean coordinates $(x^{i})_{1\leq i\leq n+1}$
\eq{  V^{2}\bar g^{ij}D_{ij}(\tfr fV) - V^{2}\bar g^{ij}\bar\Ga^{k}_{ij}D_{k}(\tfr fV)+ 2V\bar g^{ij}D_{i}VD_{j}(\tfr fV)=V,}
where $D$ is the standard connection on the Euclidean $\bbR^{n+1}$ and $(\bar\Ga^{k}_{ij})$ are the Christoffel symbols of $\bar g$ in these coordinates. On $B$, bounds and ellipticity of the coefficients $(a^{ij})$ and $(b^{i})$ are controlled by $\min_{M}r$ and $\max_{M}r$. In addition, due to the uniform equivalence of $(\bar g_{ij})$ and the Euclidean metric $(\de_{ij})$, within $B$ we can find a Euclidean ball $\ti B$, which touches $M$ at $x$. Now, working entirely in the Euclidean setting, we can quantify Hopf's boundary point lemma in a similar way as performed in the spaceform case in \cite[Lemma~4.1]{Scheuer:03/2021}. We obtain for the Euclidean normal $\ti \nu$ and using $f_{|M} = 0$
\eq{f_{\nu}\geq cf_{\ti\nu} = V(\tfr fV)_{\ti\nu}\geq c>0 ,}
with $c$ depending on $\rho^{-1}$, $\max_{M}r$ and $\min_{M}r$.

Inserting the solution $f$ into \eqref{lixia}, noticing that $V=H_1=0$ and $\nu=-\del_r$ on $\{r=0\}\times \bbS^n$, we get
\eq{\label{xeq1}
           &\frac{n}{n+1}\int_{\Om}V -\int_{\Om} V\Big|\mr{\bar \nabla}^2 f-\frac{\mr{\bar \nabla}^2 V}{V}f\Big|^2 
           \\ \ge& \int_M nV H_1f_\nu^2-c_0^2\int_{\bbS^n}V_{-\del_r}\frac{\bar \nabla^2 V(-\del_r, -\del_r)-\bar \Delta V}{V} \vt(0)^n
           \\=& \int_M nV H_1f_\nu^2-\frac{n}{(n+1)^2}\vt(0)^{1+n}|\bbS^n|.
}
On the other hand,
\eq{\label{xeq2}
           \int_{\Om}V &=\int_{\Om} V\Big(\bar \De f-\frac{\bar \De V}{V} f\Big)
= \int_{\del \Om} (V f_\nu-fV_\nu)
          \\ &= \int_M V f_\nu-c_0 \int_{\bbS^n}V_{-\del_r}\vt(0)^n=\int_M V f_\nu-\frac{1}{n+1}\vt(0)^{1+n}|\bbS^n|
}
and 
\eq{\label{xeq3}
           (n+1)\int_{\Om}V =&\int_{\Om} \ov{\dive}(\vt(r)\del_r)\\=&\int_M \bar g(\vt(r)\del_r, \nu)+\int_{\bbS^n} \bar g(\vt(0)\del_r, -\del_r)\vt(0)^n
           \\=&\int_M \bar g(\vt(r)\del_r, \nu)- \vt(0)^{1+n}|\bbS^n|.
}
Thus 
\eq{\label{xeq4}
          \int_M \bar g(\vt(r)\del_r, \nu)= (n+1)\left(\int_{\Om}V+\frac{1}{n+1}\vt(0)^{1+n}|\bbS^n|\right) =(n+1)\int_M V f_\nu. }
Define \eq{\delta=\frac{\int_M \frac{V}{H_1} }{\int_M \bar g(\vt(r)\del_r, \nu)}-1.
}
Combining \eqref{xeq1}-\eqref{xeq3} and using H\"older's inequality, we get
\eq{\label{xeq5}
           &\left( \int_M \bar g(\vt(r)\del_r, \nu)\right)^2=(n+1)^2\left(\int_M V f_\nu\right)^2
           \\ \le~& (n+1)^2\int_M \frac{V}{H_1} \int_M V H_1f_\nu^2
\\ \le~&(1+\delta)\int_M \bar g(\vt(r)\del_r, \nu)\left( \int_M \bar g(\vt(r)\del_r, \nu)- \frac{(n+1)^2}{n}\int_{\Om} V\Big|\mr{\bar \nabla}^2 f-\frac{\mr{\bar \nabla}^2 V}{V}f\Big|^2\right).
}
It follows that 
\eq{\label{xeq6'}
           \int_{\Om} V\Big|\mr{\bar \nabla}^2 f-\frac{\mr{\bar \nabla}^2 V}{V}f\Big|^2& \le\frac{n}{(n+1)^2}\frac{\delta}{1+\delta}\int_M \bar g(\vt(r)\del_r, \nu)
           \\&=\frac{n}{(n+1)^2(1+\delta)}\int_M \left(\frac{V}{H_1} -\bar g(\vt(r)\del_r, \nu)\right).
           }
              Let $r_0=\min_M r$ and $\Om_{r_0}=\{p\in \Om: r(p)\ge r_0\}$.
         Then 
         \eq{\label{xeq6}
           \int_{\Om_{\frac{r_0}{2}}} \Big|\mr{\bar \nabla}^2 f-\frac{\mr{\bar \nabla}^2 V}{V}f\Big|^2& \le \frac{n}{(n+1)^2(1+\delta)}\frac{1}{\min\limits_{r\ge \frac{r_0}{2}} V}\int_M \left(\frac{V}{H_1} -\bar g(\vt(r)\del_r, \nu)\right).
 }

Define \eq{\ep=\int_{\Om_{\frac{r_0}{2}}} \Big|\mr{\bar \nabla}^2 f-\frac{\mr{\bar \nabla}^2 V}{V}f\Big|^2.}
We may assume $\ep<1$.
As in the proof of Theorem \ref{main thm}, by using the flow \eqref{level-set-flow}, we construct a one-sided open neighbourhood of $M$ of the form 
\eq{\label{level-set-nbhd}\cU = \bigcup_{0<t<T}\{f_{|\cU}=-t\}\equiv \bigcup_{0<t<T}M_{-t},}
where  
\eq{T = \min\Big\{\fr{\min_{M}\abs{\bar\n f}^{2}}{4\|\bar\n^{2}f\|_{\8,\Om}}, \fr{\min_{M}\abs{\bar\n f}}{4}{r_{0}}, \ep^\frac{1}{1+\beta}
\Big\}.}
As calculated in \cite[p.~17/18]{Scheuer:03/2021}, on $\cU$ we have
\eq{\abs{\bar\n f}\geq \tfr12\min_{M}\abs{\bar\n f}\q\mbox{and}\q C_{n}^{-1}\abs{M_{t}}\leq \abs{M}\leq C_{n}\abs{M_{t}}.}
In addition, we need to prove that $\cU$ does not touch the other boundary component. Note that this step was not necessary in the proof of \cite[Thm.~1.1]{Scheuer:03/2021}, because there $f_{|\del\Om} = 0$. However, with our definition of $T$ and the flow \eqref{level-set-flow},
we estimate for $0<t<T$ and $\xi \in M$,
\eq{{d_{\bar M}}(F(t,\xi),\xi)\leq \int_{0}^{t}\abs{\del_{t}F(s,\xi)}\,ds\leq \tfr 12 {r_{0}}.}
Hence $\cU\sub \Om_{\frac{r_0}{2}}$.

Using the co-area formula in $\cU$, we get
\eq{\int_{-T}^{0}\int_{M_{s}}\abs{\bar\n f}^{}\abs{\mr{A}}^{2}\,ds
&\leq \int_{\{-T<t<0\}}\abs{\mr{\bar\n}^{2}f}^{2}
\\&\leq  {2}\int_{\Om_{\frac{r_0}{2}}} \Big|\mr{\bar \nabla}^2 f-\frac{\mr{\bar \nabla}^2 V}{V}f\Big|^2+{2}\int_{\{-T<t<0\}} \Big|\frac{\mr{\bar \nabla}^2 V}{V}f\Big|^2
\\&
\leq \ep+ \left\|\frac{\bar \n^2 V}{V}\right\|_{\8,\Om} T^2\int_{-T}^0\int_{M_{s}}\frac{1}{\abs{\bar\n f}}
\\&\leq \ep+CT^3.
}
Hence there exists $s\in [-T, 0]$, such that 
\eq{\|\mr{A}\|^{2}_{2,M_{s}}\leq \fr{C\ep}{T}+CT^2.}
Following the same argument as in the proof of Theorem \ref{main thm}, taking into account of the regularity estimate \eqref{reg-f} of $f$, we conclude that
\eq{\|\mr{A}\|^{2}_{2, M}\leq \fr{C\ep}{T}+CT^2+CT^\be\le C\ep^{\frac{\beta}{1+\beta}}.
}
This completes the proof of Theorem \ref{main:HK}.

\subsection*{Stability for almost CMC hypersurfaces in warped products}
We continue with quantifying Brendle's CMC classification from \cite{Brendle:06/2013}.

\begin{proof}
Let $f$ be a solution to \eqref{dirichlet} and then use $f$ in \eqref{lixia}, we get  from \eqref{xeq1} and \eqref{xeq2} that
\eq{\label{xxeq1}
           &\int_{\Om} V\Big|\mr{\bar \nabla}^2 f-\frac{\mr{\bar \nabla}^2 V}{V}f\Big|^2 
           \\ \le &\frac{n}{n+1}\left(\int_{\Om}V +\frac{1}{n+1}\vt(0)^{1+n}|\bbS^n| \right)- \int_M nV H_1f_\nu^2
           \\ =&\frac{n}{n+1}\frac{1}{\int_{\Om}V +\frac{1}{n+1}\vt(0)^{1+n}|\bbS^n|}\left(\int_M V f_\nu\right)^2- \int_M nV H_1f_\nu^2
           \\ \le&\frac{n}{n+1}\frac{\int_M V}{\int_{\Om}V +\frac{1}{n+1}\vt(0)^{1+n}|\bbS^n|} \int_M V f_\nu^2- \int_M nV H_1f_\nu^2
           \\=&n\int_M V (\mathcal{H}-H_1)f_\nu^2,
}
where 
\eq{\mathcal{H}=\frac{\int_M V}{{ (n+1)\int_{\Om}V +\vt(0)^{1+n}|\bbS^n|}}.}
It follows that 
\eq{\label{xxeq2}
           \int_{\Om} V\Big|\mr{\bar \nabla}^2 f-\frac{\mr{\bar \nabla}^2 V}{V}f\Big|^2 \le C\left\|H_1-\mathcal{H}\right\|_{1,M}.
}
Arguing as in the proof of Theorem \ref{main:HK}, we get
\eq{\label{xxeq4}\|\mr{A}\|^{2}_{2, M}\leq C\left(\int_{\Om} V\Big|\mr{\bar \nabla}^2 f-\frac{\mr{\bar \nabla}^2 V}{V}f\Big|^2 \right)^{\frac{\beta}{1+\beta}}\le C\left\|H_1-\mathcal{H}\right\|_{1,M}^{\frac{\beta}{1+\beta}}.}

Next, by the Codazzi equation, we have
\eq{
\n^{j}\mathring{A}_{ij}=\frac{n-1}{n}\n_i H+\bar R_{i\nu},
}
where $\bar R_{i\nu}=\ov{\Rc}(e_i,\nu)$.
By using the divergence theorem, we get 
\eq{
\int_M\sum_i \bar R_{i\nu}^2&=-\int_M \br{(n-1)\n^{i} H_{1} \bar R_{i\nu}- (\n^{j}\mathring{A}_{j}^i)\bar R_{i\nu}} 
\\&=\int_M \br{(n-1) (H_1-\mathcal{H}) \n^{i}\bar R_{i\nu}-\mathring{A}_{j}^{i}\n^{j}\bar R_{i\nu}}
\\&=\int_M (n-1)  (H_1-\mathcal{H}) (\bar \n^{i}\bar R_{i\nu}- H\bar R_{\nu\nu}+h^{ij}\bar R_{ij})\\
&\hp{=}-\mathring{A}_{j}^{i}(\bar \n^{j}\bar R_{i\nu}- h^{j}_{i}\bar R_{\nu\nu}+h_{k}^{j}\bar R^{k}_{i})
\\&\le C(\|H_1-\mathcal{H}\|_{1,M}+\|\mathring{A}\|_{1,M}).
}

By using the explicit formula \eqref{Rc} for $\ov{\Rc}$, the assumption \eqref{H4} and the estimate \eqref{xxeq4},  we deduce
\eq{\label{xxeq5}
\int_M |\del_r^T|^2\left<\del_r,\nu\right>^2\le C(\|H_1-\mathcal{H}\|_{1,M}+\|\mathring{A}\|_{1,M})\le C\|H_1-\mathcal{H}\|_{1,M}^{\frac{\beta}{1+\beta}},
}
where we assume again without loss of generality that $\|H_{1}-\cH\|_{1,M}<1$.

Now we claim that there exists some $\delta_0>0$,  depending on $|M|_{C^2}$, such that if $\|H_1-\mathcal{H}\|_{1,M}\le \delta_0$, then \eq{\label{starshaped}\left<\del_r,\nu\right> \ge \frac14\hbox{ on }M.}
If the claim is not true, then there exists some point $p\in M$ such that $\left<\del_r,\nu\right>(p) < \frac14$. 

At a maximal point of $r_{|M}$ we have $\nu = \del_{r}$. Thus there exists at least one point in $M$ at which $\ip{\del_{r}}{\nu}=1$. By continuity, there exists some other point $p_1$ such that \eq{\left<\del_r,\nu\right>(p_1)=\tfr12.}Hence
\eq{|\del_r^T|^2\left<\del_r,\nu\right>^2(p_1)=(1-\left<\del_r,\nu\right>^2)\left<\del_r,\nu\right>^2 (p_1)=\tfr{3}{16}.
}
It follows that there exists some $\delta_1>0$, depending on $|M|_{C^2}$, such that  for all $q\in M $ with $d_{M}(q,p_{1})<\delta_{1}$ there holds
\eq{|\del_r^T|^2\left<\del_r,\nu\right>^2(q)\ge \tfrac{1}{16}.}
This would contradict \eqref{xxeq5} for $\|H_1-\mathcal{H}\|_{1,M}$ small enough.
Therefore, we conclude from \eqref{starshaped} that when $\|H_1-\mathcal{H}\|_{1,M}\le \delta_0$, 
\eq{\label{xxeq6}
\int_M (1-\left<\del_r,\nu\right>^2)=\int_M |\del_r^T|^2 \le 16\int_M |\del_r^T|^2\left<\del_r,\nu\right>^2\le C\|H_1-\mathcal{H}\|_{1,M}^{\frac{\beta}{1+\beta}}.
}
Finally, since $\left<\del_r,\nu\right>>0$ on $M$, we can represent $M$ as a graph over $\{r=0\}\times \bbS^n$ with graph function $\varphi$. Then
\eq{\label{xxeq7}
1-\left<\del_r,\nu\right>^2=\frac{ \vt^{-2}|\n^\si \varphi|^2}{1+\vt^{-2}|\n^\si \varphi|^2},
}
where $\n^\si$ denotes the covariant derivative on $(\bbS^n, \si)$.
By the Poincare inequality on $\bbS^n$, we deduce from \eqref{xxeq6} and \eqref{xxeq7} that
\eq{
\int_{\bbS^n} |\varphi-\bar \varphi|^2\le \frac1n\int_{\bbS^n} |\n^\si  \varphi|^2\le C\|H_1-\mathcal{H}\|_{1,M}^{\frac{\beta}{1+\beta}}.
}
where $\bar \varphi=\frac{1}{|\bbS^n|}\int_{\bbS^n}\varphi$.

Note that $M$ is of $C^{2, \beta_2}$. By the Gagliardo-Nirenberg interpolation inequality and the Sobolev-Morrey embedding theorem, we get immediately for any $p>n$ and $\gamma=1-\frac{n}{p}$, there exists some $\alpha(n)$, depending only on $n$, such that %
\eq{\label{pf:thm-CMC-1}
\abs{\varphi-\bar \varphi}_{0,\ga,\bbS^{n}}\le C\|\varphi-\bar \varphi\|_{W^{1,p}}\le C\|H_1-\mathcal{H}\|_{1,M}^{\alpha(n)\beta}.
}
This clearly implies the assertion.
\end{proof}

\providecommand{\bysame}{\leavevmode\hbox to3em{\hrulefill}\thinspace}
\providecommand{\MR}{\relax\ifhmode\unskip\space\fi MR }
\providecommand{\MRhref}[2]{%
  \href{http://www.ams.org/mathscinet-getitem?mr=#1}{#2}
}
\providecommand{\href}[2]{#2}


\end{document}